\documentclass[a4paper,12pt]{article}
\usepackage{amsfonts}
\usepackage{amssymb}
\usepackage{latexsym}
\usepackage{amsmath}
\usepackage{amscd}
\usepackage{amsthm}
\usepackage{graphicx}
\usepackage{epsfig}
\usepackage{hyperref}
\usepackage{url}
\usepackage{xcolor}

\newtheorem{theorem}{Theorem}[section]

\newtheorem{proposition}[theorem]{Proposition}

\newtheorem{remark}[theorem]{Remark}

\newcommand{\N}{{\mathbb N}}
\newcommand{\Z}{{\mathbb Z}}
\renewcommand{\P}{{\mathbb P}}

\newcommand{\R}{{\mathbb R}}

\newcommand{\x}{{\mathbf x}}
\newcommand{\y}{{\mathbf y}}

\newcommand{\w}{{\mathbf w}}

\title{ Co-existence in a discrete time competing Frog Model}
\author{Rishideep Roy \& Kumarjit Saha }

\begin{document}
	
	\maketitle
	
	\begin{abstract}
		
		We study coexistence in discrete time multi-type frog models.
		We first show that for two types of particles on $\mathbb{Z}^d$, for $d\geq2$, for any jumping parameters $p_1, p_2 \in (0,1]$, coexistence occurs with positive probability for sufficiently rich deterministic initial configuration. We extend this to the case of random distribution of initial particles. We study the question of coexistence for multiple types and show positive probability coexistence of $2^d$ types on $\mathbb{Z}^d$ for rich enough initial 
		configuration. We also show an instance of infinite coexistence on $\mathbb{Z}^d$ for $d \geq 3$ provided we have sufficiently rich initial configuration.
	\end{abstract}

\section{Introduction and model}
Coexistence in frog models, in the study of interacting particle systems has been an object of recent interest. The frog model can be viewed as a model for describing information
spreading. The original idea is that every active particle has some information and it shares that information with a sleeping particle at the time the former meets the latter. Particles that have the information move freely
helping in the process of spreading information. The name `frog model' seems to be attributed to R. Durrett in the literature. The question of recurrence for frog model on $\mathbb{Z}^d$ was studied in \cite{TW99} and \cite{P01}. The question of extinction and survival was studied in \cite{AMP02} for a modification of the frog model, where  active particles may disappear at each step. 
The shape theorem for the frog model on $\mathbb{Z}^d$ was studied in \cite{AMP02a}, \cite{AMPR01} and extended for a continuous time version of the model in \cite{RS04}.

Coexistence in discrete time two-type competing frog model has been recently studied 
in \cite{DHL19} with an assumption of equality of the two jump probabilities. 
In this paper we study the problem of coexistence for two types 
without the equality assumption, as well as coexistence for multi types and 
show that,  in the beginning if each of the dormant sites has sufficient number of particles
with high probability then the coexistence probability is positive.  
It may be worthy to note that the work of Deijfen et.\ al in \cite{DHL19} has been extended to show coexistence for any values of the jump probabilities for dimension $1$ in \cite{HK20}.

The two-type competing frog model was introduced in \cite{DHL19}.
Let $\{\eta(\x) : \x \in \Z^d \}$ denote a family of i.i.d. 
non-negative integer valued random variables
such that $\eta(\x)$ denotes the initial number of particles (frogs) at site $\x$ and 
these particles (frogs) are in the dormant state. 
At time $0$, particle at the origin is activated and assigned type $1$,
while particle at another site $\y \in \Z^d \setminus \{\mathbf{0}\}$ 
is activated and assigned type $2$. We will call a site $\x \in \Z^d \setminus \{\mathbf{0}, \y\}$ as initially dormant site. 
Fix $p_1, p_2 \in (0,1]$ and active particles of both types follow random walks as follows.  
At each time point,  for $i=1,2$, each of the $i$-type active particles 
independently either stay back at their present site with probability $1 - p_i$, 
or move to a different site,  uniformly chosen among all possible neighbouring sites,
with probability $p_i$. The active particles on reaching any site 
activate all the dormant particles present there and assign them the type of itself.
Clearly, an active type $i$ particle follows simple symmetric {\it lazy} 
random walk if $p_i < 1$. Keeping this in mind, from hereafter with a slight abuse of notation,
we call $p_1, p_2 \in (0,1]$ as laziness parameters of the two-type competing frog model 
described above. The resulting random walks followed by active particles 
will be called {\it possibly lazy} nearest neighbour random walks. 

We now describe a construction of the two-type competing frog model
along the lines of Deijfen et. al. \cite{DHL19}. 
Let $\{S^{\x, j, (i)} : \x \in \Z^d, j \geq 1 \}_{i = 1, 2}$ denote 
independent collections of nearest neighbour possibly lazy 
random walks on $\Z^d$, starting from the origin, such that 
the distribution of the laziness clock is given by:
$$
\P(\min \{ n \geq 1: S^{\x, j, (i)}_n \neq S^{\x, j, (i)}_0 \} = k) = (1 - p_i)^{k-1}p_i \text{ for } k \in \N.
$$
We further assume that both the collections of random walks are independent of 
the initial configuration $\{\eta(\x) : \x \in \Z^d\}$.  
%In what follows, we will call both $p_1$ and $p_2$ as the `laziness' parameters of our model. 
%Using the above defined independent collections of random walks, 
%movement of active particles are modelled as follows: 

We say that a site is discovered when it is first hit by an active particle.
Suppose at time $n$, a dormant site $\x \in \Z^d$ is discovered by type `$i$'
active particle(s) only. 
Then the $j$-th  `newly' activated particle at $\x$, which has to be of type '$i$',
follows random walk $\{S^{\x, j, (i)}_{n+l}\}_{l \geq 1}$, and is at the site
$(\x + S^{\x, j, (i)}_{n+l})$ at time $n+l$, for $l \geq 1$. In case $\x$ has been discovered 
by both types of active particles at the same instant, some tie breaking mechanism will 
be used to decide the type of the newly activated particles.  
We are not specifying any tie-breaking mechanism here as our 
result remains valid for any tie-breaking mechanism.

We assume that at time $0$ we start with one active type $1$ particle at the origin, and one active type $2$ particle at some other site $\y $, and the corresponding conditional probability measure is denoted by $\P_{\mathbf{0}, \y}(.)$. 
For $i = 1,2 $, let $N_n(i)$ denote the set of $i$-type active particles at time $n$ and let $|N_n(i)|$ denote the cardinality of the same. By our assumption, we have 
$\P_{\mathbf{0}, \y}(|N_0(1)|= |N_0(2)| = 1)$. 
For $i = 1, 2$ let $G_i$ denote the event that 
$$
G_i := \{ \mid N_n(i) \mid \to \infty \text{ as }n \to \infty \}.
$$
For two-type competing frog model the `coexistence' event 
is expressed as the event $G_1\cap G_2 $. Under certain assumptions on the initial configuration $\{\eta(\x) : \x \in \Z^d\}$, 
in \cite{DHL19} it is proved that for $p_1 = p_2$, coexistence occurs with positive probability. 
A natural question of interest is what happens when $p_1 \neq p_2$. It has been conjectured in \cite{DHL19} that if $\eta(\x)$ is a heavy-tailed random variable,
then coexistence is possible. The reason behind this intuition is that,
irrespective of the values of $p_1$ and $p_2$, 
both types will have the same limiting shape: a full diamond $D := \{\x \in \R^d : ||\x||_1 \leq 1\}$
(Theorem 1.5 of \cite{DHL19}). 

In this paper in Theorem \ref{thm:CoexistsenceHighProb}, 
we show positive probability of coexistence for any $p_1 , p_2 \in (0,1]$
and  for a deterministic initial configuration provided each dormant site 
has sufficiently many  particles in the beginning. 
We neither require equality of the entire limiting shapes nor heavy tailed 
initial configurations. In fact, we prove coexistence for
bounded (but sufficiently large) initial configurations. 
The main tool for our proof is a coupling with oriented 
(site) percolation on $\Z^d$, defined in Section \ref{sec:CouplingOP}. 
In Proposition \ref{prop:CoexistencewithZeroParticle} we extend this for random i.i.d. 
initial configuration. The probability of having sufficient number of particles 
per site must be large though.
We mention here that for any $p_1 = p_2 $,
positive probability of coexistence was shown in \cite{DHL19} for one particle per site 
initial configuration as well. Our results are not applicable for one particle per site 
initial configuration.

We apply this method to study coexistence for more than two types of particles as well. 
Coexistence of more than two types is denoted by the event that
the number of active particles of each type
grows to infinity simultaneously. We obtain positive 
probability for coexistence of $2^d$ types on $\Z^d$ for $d \geq 2$ and of 
{\it infinitely} many types of particles on $\mathbb{Z}^d$ for $d \geqslant 3$.
These results are applicable for any laziness parameters provided 
initially, each dormant site has sufficient number of particles with high probability. 

The succeeding part of our paper is divided into three sections. The first of them, Section \ref{sec:CouplingOP} defines a coupling of (possibly) lazy frog model 
with an oriented (site) percolation model. 
This coupling is the basis of all our proofs.  
In Section \ref{sec:Coex_Gen_p_Multitype_1} we prove Theorem \ref{thm:CoexistsenceHighProb} and 
Proposition \ref{prop:CoexistencewithZeroParticle}. 
In Section \ref{sec:Coex_Multi-types} we 
study the question of coexistence for general multi-type frog models and show positive probability of coexistence of $2^d$ types on $\mathbb{Z}^d$. We also show  coexistence of 
infinitely many types for $d \geqslant 3$. 

\section{A coupling with oriented percolation}
\label{sec:CouplingOP}

In this section we define a coupling of the (possibly lazy) frog model 
with an oriented site percolation model on $\mathbb{Z}^d$. 

To this end, we consider collection of i.i.d. random tuples  
$\{ (I^{\x, j}_{n}, U^{\x, j}_{n}) : \x \in \Z^d, j \geq 1\}_{n \geq 0}$, independent 
of the initial i.i.d. configuration $\{\eta(\x) : \x \in \Z^d \}$, 
such that $U^{\mathbf{0}, 1}_{0}$ is $U(0,1)$ random variable and $I^{\mathbf{0}, 1}_{0}$ denotes the increment for a nearest neighbour random work on $\Z^d$, i.e.,
$$
\P(I^{\mathbf{0}, 1}_{0} = + e_i) = \P(I^{\mathbf{0}, 1}_{0} = - e_i)  = 1/2d \text{ for all } 1\leq i \leq d,
$$  
where $\{e_1, \cdots, e_d \}$ denotes the standard basis set for $\R^d$. We further assume that the collections, $\{ I^{\x, j}_{n} : \x \in \Z^d, j \geq 1\}_{n \geq 0}$ and $\{ U^{\x, j}_{n} : \x \in \Z^d, j \geq 1\}_{n \geq 0}$, are mutually independent. 

Fix $p \in (0,1]$ and we will construct a frog model with laziness parameter $p$
using the collection $\{ (I^{\x, j}_{n}, U^{\x, j}_{n}) : \x \in \Z^d, j \geq 1\}_{n \geq 0}$.  
At time $t$, consider an active site $\x$ (which has been activated at some earlier 
time $0 \leq s \leq t$). 
We further assume that we have an algorithm which orders the {\it active} particles present at $\x$ 
at time $t$. Note that, some of these active particles may come from some other sites as well.  
The $j$-th active particle present at site $\x$ at time $t$ 
jumps to the site $\x + I^{\x, j}_{t}$ at time $t + 1$ only if
$U^{\x, j}_{t} \leq p$, otherwise it stays still. 
We observe that on the event $\{U^{\x, j}_{t} \leq p \}$, the $j$-th active particle at $\x$ at time $t$ reaches site $(\x + I^{\x, j}_{t})$ at time $t+1$ and it's movement (if any) for the next time point is decided by the collection of random vectors 
$\{(U^{(\x + I^{\x, j}_{t }), j^\prime}_{t + 1}, I^{(\x + I^{\x, j}_{t}), j^\prime}_{t + 1}) : j^\prime \geq 1 \}$. 

This describes a frog model with laziness parameter $p$. Further, the use of $U(0,1)$
random variables allows us to couple frog models with different laziness parameters. 
Consider a two-type competing frog model where type $i$ has laziness parameter $p_i$
for $i = 1, 2$. Then the $j$-th active particle  of type $i$ present at site $\x$ at time $t$
jumps to the site $\x + I^{\x, j}_{t}$ at time $t + 1$ only if
$U^{\x, j}_{t} \leq p_i$, otherwise it stays still.  
This describes the same process as considered in \cite{DHL19}. Heuristically, in \cite{DHL19}, to each newly activated particle an independent possibly lazy random walk trajectory was attached, depending upon its type, whereas in our construction each active particle, upon reaching a new location, gets an i.i.d. (possibly) lazy increment, depending upon its type.

We now use the above construction of possibly lazy frog model to 
couple it with an oriented site percolation process.  
For $\x \in \Z^d$ and $\theta \in \{+1, -1\}^d$ we define the  `$\theta$' orthant 
starting from $\x$ as 
\begin{align*}
\Lambda^\theta(\x) := \{\x + \sum_{j=1}^d k_j\theta(j) e_j  :  k_j \in \N \cup\{0\} 
\text{ for all }1 \leq i \leq d \},
\end{align*}
where $\theta(j)$ denotes the $j^{th}$ co-ordinate of $\theta$.
By definition we have $\x \in \Lambda(\x)$. 
Fix $M \in \N$ and $\ell \in \N \cup \{0\}$. Corresponding to the choice of $M, \ell$ and 
$p \in (0,1]$ we define an oriented (site) percolation model on the 
orthant $\Lambda^\theta(\x)$.
For a site $\w $ in the orthant $\Lambda^\theta(\x)$ with 
$||\x - \w ||_1 = a $, we say that  $\w$ is `open' if
the following event occurs
\begin{align*}
\{ \w \text{ is  open}\} := & \{\eta(\w) \geq M\} \bigcap  \\
& \qquad \Bigl( \cap_{j=1}^d  
\{ I^{\w,i}_{\ell + a} =  \theta(j)e_j \text{ with }U^{\w,i}_{\ell + a} \leq p \text{ for some }1\leq i \leq M+1 \} \Bigr ).
\end{align*}
In other words, for an open site $\w \in \Lambda^\theta(\x) $ with 
$||\x - \w||_1 = a$, 
if there are at least $M+1$ {\it active} particles at $\w$ at time $\ell + a$, then 
all the oriented neighbours of $\w$, i.e., $\w + \theta(j)e_j$ for all $1 \leq j \leq d$ are reached 
by some of these  active $M+1$ particles at time $\ell + a+ 1$. 
%It is important to observe that contrary to the usual oriented percolation model, for $\x \in \Lambda(\theta_m)$ the notion of oriented neighbours of $\x$ depends on the orthant $\Lambda(\theta_m)$
%to which it belongs to and given by $\{\x + \theta(j)e_j 
%: 1 \leq j \leq d\}$.   

Let $\mathbf{1}_{\{\w \text{ open}\}}$ denote 
the indicator random variable 
corresponding to the event that $\w$ is open. 
We observe that the collection $\{\mathbf{1}_{\{\w \text{ open}\}} : 
\w \in \Lambda^\theta(\x)\}$ gives an i.i.d. collection of Bernoulli random variables 
with success probability given by $\P(\eta(\w ) \geq M )g(M, p)$ where $g(M,p)$
is given by
\begin{eqnarray*}g(M, p) &=& \sum_{\substack{c_1,c_2,\ldots,c_d \geq 1, \\\sum_{i=1}^d c_i \leq M+1} }
	\dfrac{(M+1)!}{c_1! c_2! \ldots c_d !(M+1-c_1 - c_2 -\cdots -c_d)!} \\
	&& \qquad \times  \left( \frac{p}{2d}\right)^{c_1+c_2+\cdots+c_d}\left( 1 -\frac{p}{2} 
	\right)^{ M + 1 - c_1 - c_2 -\cdots -c_d} \\
	& \geq & 1 - d( 1 - \frac{p}{2d})^{ M + 1 }.
\end{eqnarray*}
We end this section with the observation that for any $p \in (0,1]$, we can choose  $M = M(p)$ large 
to make $g(M, p)$ arbitrarily close  to $1$.

\section{Coexistence for two types}
\label{sec:Coex_Gen_p_Multitype_1}

In this section we prove coexistence for two-type competing frog model for any $p_1, p_2 \in (0,1]$.  
We first prove this for a deterministic initial configuration such that  
each dormant site has sufficiently large number of particles. 
We later extend this result for random i.i.d. initial configuration. 

%Below we state our first result which says that for any $p_1 , p_2 \in (0,1)$ then there exists $M = M(p_1 \wedge p_2, d)$,
%which depends only on $p_1 \wedge p_2$ and $d$, 
%so that for initial configuration $\eta(\x) \geq M$ for all $\x \in \Z^d \setminus \{\mathbf{0}, \y\}$, 
%coexistence occurs with positive probability.
We recall that at time $0$, we start the two-type competing frog model 
with only two active particles, one type $1$ active particle
at the origin and another type $2$ active particle at some other site $\y $.
The corresponding conditional probability measure is denoted by $\P_{\mathbf{0}, \y}(\cdot)$. The following is our first result in this section: 
\begin{theorem}
	\label{thm:CoexistsenceHighProb}
	For any $p_1, p_2 \in (0,1]$, there exists $M = M(p_1 \wedge p_2, d)$ such that  for  
	initial configuration $\eta(x) \geq M$ for all $x \in \mathbb{Z}^d \setminus \{\mathbf{0} , \y\}$ deterministically, 
	coexistence occurs with positive probability. 
\end{theorem}

\noindent \textbf{ Proof:} Our proof is based on the coupling defined in Section 
\ref{sec:CouplingOP}. Let us call $\underset{\sim}{e}=(1,1,\ldots,1)$. We choose $m \in \N$ such that 
$||\y ||_1 \leq m$ and we consider two vertices $\theta^+_m, \theta^-_m \in \Z^d$, where, 
$$
\theta^+_m := (+m , \cdots , +m) \text{ and } \theta^-_m := (-m , \cdots , -m).
$$
We consider two orthants given as
\begin{align*}
\Lambda^{\underset{\sim}{e}}(\theta^+_m) & := \{ \theta^+_m + \sum_{j=1}^d k_je_j  :  k_j \in \N \cup\{0\} 
\text{ for all }1 \leq i \leq d \} \text{ and }\\
\Lambda^{-\underset{\sim}{e}}(\theta^-_m) & := \{ \theta^-_m - \sum_{j=1}^d k_je_j  :  k_j \in \N \cup\{0\} 
\text{ for all }1 \leq i \leq d \}.
\end{align*} 
By an abuse of notation we shall call them $\Lambda(\theta^+_m)$ and $ \Lambda(\theta^-_m)$ respectively. Clearly, the two orthants  $ \Lambda(\theta^+_m)$ and $ \Lambda(\theta^-_m)$
are disjoint.  

A finite sequence of nearest neighbour lattice points 
$\{\x_i\}_{0 \leq i \leq n} \subseteq \Z^d$ gives us a `path' of length $n$ starting from 
$\x_0$ and ending at $\x_n$. 
Let $\pi_1 := \{\x_0 := \mathbf{0}, \x_1, \cdots , \x_{md -1}, \x_{md} := \theta^+_m\}$ 
be a path of length $md$ from $\mathbf{0}$ to $\theta^+_m$ such that $||\x_i||_1 < md $
for all $0 \leq i \leq md - 1$. Similarly, set $k = ||\y - \theta^-_m||_1$ and 
let $\pi_2 := \{\y_0 := \mathbf{y}, \y_1, \cdots , \y_{k -1},\y_k := \theta^-_m\}$ 
denote a path of length $k$ from $\mathbf{y}$ to $\theta^-_m$  such that $||\y_i||_1 < md $
for all $0 \leq i \leq k - 1$. Clearly, such $\pi_1$ and $\pi_2$ exist. 

We choose  $M(p_1\wedge p_2, d) \in \N$ large such that $g(M, p_1\wedge p_2) > p^\uparrow_c $, where
$p^\uparrow_c = p^\uparrow_c(d)$ is the oriented site percolation threshold for $\mathbb{Z}^d$.
We use the coupling defined in Section \ref{sec:CouplingOP} for the 
orthant $\Lambda(\theta^+_m)$ and consider the corresponding oriented 
percolation model in $\Lambda(\theta^+_m)$ for the choice 
$p = p_1\wedge p_2, M = M(p_1\wedge p_2, d)$ and 
$\ell = md$. On the other hand, for the other orthant $\Lambda(\theta^-_m)$ 
we consider the corresponding oriented 
percolation model with $p = p_1\wedge p_2, M = M(p_1\wedge p_2, d)$ and 
$\ell = k$. We assume that in the beginning, any site $\x \in \Z^d 
\setminus \{\mathbf{0}, \mathbf{y}\}$ has at least $M$ many particles,
i.e., $\P(\eta(\x) \geq M) = 1$ for all $\x \in \Z^d 
\setminus \{\mathbf{0}, \mathbf{y}\}$. 

For $i = 1,2$ and $n \geq 0$, let $L_n(i) \subset \Z^d$ denote
the collection of location(s) of all $i$-type active particles at time $n$. 
Clearly, we have $\P_{\mathbf{0}, \y}(L_0(1) = \{\mathbf{0}\}, L_0(2) = \{\mathbf{y}\}) = 1$.    
Now we are ready to define the following events:
\begin{align*}
A_1 & := \{ \text{The initial type 1 active particle at the origin 
	follows path }\pi_1 \text{ and reaches } \theta^+_m \\
& \qquad \text{ at time } md \text{ and  
	the initial type 2 active particle at }\mathbf{y}
\text{ follows the path } \pi_2 \text{ and }\\
& \qquad \text{ reaches }\theta^-_m 
\text{ exactly at time } k\}.   \\
A_2 & := (\cap_{n = 0}^{md} L_n(1)  = \{\x_n \}) \bigcap (\cap_{n = 0}^{k} L_n(2)  = \{\y_n \}).\\ 
A_3 & := \{ \theta^+_m \text{ is a  percolating point, i.e., 
	there exists an infinite `oriented' open path }\\
& \qquad  \text{  consisting of sequence of oriented open neighbours starting from }\theta^+_m
\text{ in }\Lambda(\theta^+_m) \}.\\
A_4 & := \{ \theta^-_m \text{ is a  percolating point} \}.
\end{align*}
We observe that the event $A_1\cap A_2$ ensures that till time $md$, 
active type $1$ particles are not allowed to be lazy and they are {\it all} 
assembled at $\theta^+_m$ at time $md $. Similarly, till time $k$  all the active type $2$
particles move without being lazy and reach the site $\theta^-_m$ at time $k$.

Set $M^\prime := \max\{\eta(\w) : \w \in \pi_1 \cup \pi_2\}$ and we observe that
on the event $A_1\cap A_2$, the number of newly activated particles 
at a site $\w \in \pi_1 \cup \pi_2$ is bounded by $M^\prime \times (k \vee md)$.
Hence, the event $A_1\cap A_2$ depends on the finite collection of random vectors 
\begin{align}
\label{def:A^1_3_Collection}
\bigl \{ (I^{\x_i, j}_{||\x_i||_1}, U^{\x_i, j}_{||\x_i||_1}), (I^{\y_j, j}_{||\y_j - \y||_1}, 
U^{\y_j, j}_{||\y_j - \y||_1}) : & 1 \leq j \leq M^\prime \times (k \vee md) + 1, \nonumber \\
& 0 \leq i  \leq md-1, 0 \leq j \leq k-1 \bigr \}.
\end{align}
Clearly, we have $\P_{\mathbf{0}, \y} ( A_1\cap A_2 ) > 0$. 
Next, we observe that for $\w \in \Lambda(\theta^+_m)$ 
we have $||\w||_1 = md + ||\w - \theta^+_m||_1$ and 
occurrence of the event $A_3$ depends on the collection 
\begin{align*}
\{(I^{\x,j}_{||\x||_1}, U^{\x,j}_{||\x||_1}) :  & \x \in \Lambda(\theta^+_m), 1 \leq j \leq M + 1\} \\
& \qquad =  \{(I^{\x,j}_{md + ||\x - \theta^+_m||_1}, U^{\x,j}_{md + ||\x - \theta^+_m||_1}) 
:  \x \in \Lambda(\theta^+_m),  1 \leq j \leq M + 1\}
\end{align*}
which is disjoint from the collection considered in (\ref{def:A^1_3_Collection}). 
Hence, the event $A_3$ is independent of the event $A_1\cap A_2$. 
By the same argument, the event $A_4$, 
which depends on the collection 
$$
\{(I^{\w,j}_{k + ||\w - \theta^-_m||_1}, U^{\w,j}_{ k + ||\w - \theta^-_m||_1}) :
\w \in \Lambda(\theta^-_m),  1 \leq j \leq M + 1\}, 
$$ 
is independent of $A_1\cap A_2$. 
The choice of $M$ as mentioned earlier ensures that 
$$
\P_{\mathbf{0}, \mathbf{y}}(A_3) = \P_{\mathbf{0}, \mathbf{y}}(A_4) > 0.
$$
Further, as the orthants $\Lambda(\theta^+_m)$ and $\Lambda(\theta^-_m)$ are disjoint, 
the events $A_3$ and $A_4$ are independent. This gives us that 
$$
\P_{\mathbf{0}, \mathbf{y}}(\cap_{i=1}^4 A_i) > 0.
$$
Finally, we show that occurrence of the event $\cap_{i=1}^4 A_i$ implies coexistence  of both types.
Event $A_1\cap A_2$ ensures that at time $md$ all the  
type $1$ active particles are at the site $\theta^+_m$ and it has at least $M + 1$
type $1$ particles. On the other hand, 
at time $k$ all the type $2$ active particles are at $\theta^-_m$
and there are at least $M+1$ type $2$ particles.  
We now claim that on the event $\cap_{i=1}^2 A_i$, a site $\x \in \Lambda(\theta^+_m)$
with $||\x - \theta^+_m||_1 = l$ {\it cannot} be reached 
by a type $2$ active particle in time $md+ l$. 
We observe that the choice of $\pi_2$ ensure that before time $k$, 
all the type $2$ active particles are inside the set $\{\w \in \Z^d : ||\w||_1 \leq md - 1\}$
and no site in $\Lambda(\theta^+_m)$ has been visited by an active
type $2$ particle. Since, all the type $2$ particles 
are at $\theta^-_m$ at time $k$ and starting from $\theta^-_m$, 
a type $2$ particle requires at least $2md + l$ many additional jumps to reach 
$\x$. Hence, an active type $2$ particle cannot reach $\x$
by time $k + 2md + l - 1$ which is strictly bigger than $md + l$. 
This proves our claim.

Similarly, a site $\w \in \Lambda(\theta^-_m)$
with $||\w - \theta^-_m||_1 = l$ {\it cannot} be reached 
by an active type $1$ particle by time $k + l$. 
The event $A_3 \cap A_4$ ensures 
that both the sites $\theta^+_m$ and $\theta^-_m$ are percolating 
points and have infinite oriented open paths 
in the orthants $\Lambda(\theta^+_m)$ and $\Lambda(\theta^-_m)$ respectively. 
We will show that dormant particles at a site $\x \in \Lambda(\theta^+_m)$ 
with $||\x - \theta^+_m ||_1 = l$,  
connected to $\theta^+_m$ through an oriented open path, 
gets activated by a type $1$ particle exactly at time $md + l$. 
We will prove this using method of induction. 

For $l = 0$, this is guaranteed by the event $A_1 \cap A_2$. 
Assuming that this is true for $l=l_0$, we show that this holds for $l= l_0 + 1$ as well. 
Consider $\x \in \Lambda(\theta^+_m)$ with $||\x - \theta^+_m||_1 = l_0 + 1$, connected through an oriented open path to $\theta^+_m$. Let $\w \in \Lambda(\theta^+_m)$ be such that:
\begin{itemize}
	\item[(i)] $\w$ is connected to $\theta^+_m$ through an oriented open path 
	(which means $\w$ must  be open) and 
	\item[(ii)] $||\w - \theta^+_m||_1 = l_0$ and $||\w - \x||_1 = 1$. 
\end{itemize}
The last condition suggests that $\x$ must be an oriented neighbour of $\w$. 
Since, $\x$ is connected to $\theta^+_m$ through an oriented open path, 
such a $\w$ must exist. Further, by our induction hypothesis, dormant particles at 
$\w$ must be activated at time $l_0$ by type $1$ particle(s) only. 
This implies that,  at time $md + l_0$ all the dormant particles at $\w$ become 
type $1$ and there are at least $M+1$ many type $1$ active particles 
present at $\w$ at time $md + l_0$.  
As $\w$ itself is an open vertex and $\x$ is an oriented neighbour of $\w$,
site $\x$ must be reached by a type $1$ active particle at time $md + l_0 + 1$
jumping from $\w$. Earlier, we proved that no type $2$ particle can reach 
$\x$ by time $md + l_0 + 1$. This implies that the dormant particles at 
$\x$ must be activated by a type $1$ particle only at time $md + l_0 + 1$. 
This completes our induction argument.
This also shows that for an infinite oriented open path starting 
from $\theta^+_m$ in $\Lambda(\theta^+_m)$, dormant particles at 
each (open) site on such a path must belong to $N_n(1)$ for some $n$. 
This implies occurrence of the event $G_1$.  

By the same argument,  there will also be an infinite \textit{oriented} 
path of open sites starting from the point $\theta^-_m$  in $\Lambda(\theta^-_m)$ and dormant particles at 
each site on such an infinite path must belong to $N_n(2)$ for some $n$. 
Hence, the event $G_2$ occurs as well. This completes the proof.
\qed

As we commented earlier, we note that Theorem \ref{thm:CoexistsenceHighProb} holds for 
{\it any}  tie-breaking mechanism. We extend the above construction 
for random i.i.d. initial configuration naturally. For i.i.d. initial configuration 
if probability of having large number of initial particles is high enough, 
we can still apply the same argument as in Theorem \ref{thm:CoexistsenceHighProb}. 

\begin{proposition}
	\label{prop:CoexistencewithZeroParticle}
	Fix $p_1, p_2 \in (0,1]$ and consider an initial configuration of i.i.d. non-negative integer valued random variables $\{\eta(\x) : \x \in \Z^d\}$.
	There exist $M=M(p_1 \wedge p_2, d)$ as in Theorem \ref{thm:CoexistsenceHighProb} and 
	$\theta \in (0,1)$, depending on $M$,  such that if $\P(\eta(\mathbf{0}) \geq M) \geq \theta$, 
	then under the probability measure $\P_{\mathbf{0}, \y}(.)$, coexistence occurs with positive probability.
\end{proposition}
%We comment here that in order to show positive probability of coexistence (for the choice $p_1 = p_2$) Deijfen et. al. \cite{DHL19} require the assumption that 
%$\eta(\x) \geq 1$ for all $\x \in \Z^d$ a.s. or $\E(\eta(\x)) < \infty$. Proposition \ref{prop:CoexistencewithZeroParticle} shows that we don't have such restrictions. On the other hand, for the initial configuration one particle per site, 
%Deijfen et. al. \cite{DHL19} was able to prove positive probability of coexistence for the choice $p_1 = p_2$ whereas our argument require sufficiently rich initial configuration.  

\noindent \textbf{ Proof:} The proof follows from the same argument as in Theorem \ref{thm:CoexistsenceHighProb} with the observation that as $g(M, p_1\wedge p_2) > p^\uparrow_c$, 
we can choose  
$\theta = \P(\eta(\x) \geq M ) \in (0,1)$ so that we have 
$\theta g(M, p_1\wedge p_2) > p^\uparrow_c $.
\qed

\begin{remark}
	\label{rem:EqualityLimitingShape}
	In \cite{DHL19}, regarding coexistence of two types of particles with different 
	$p_1$ and $p_2$, it has been conjectured that equality of the (complete) limiting 
	shapes would imply coexistence (Conjecture 1.6 of \cite{DHL19}). Extending this further, 
	the authors of \cite{DHL19} posed an interesting question: 
	whether equality of the limiting shapes along some specific direction 
	only would imply coexistence. 
	In the set up of Theorem \ref{thm:CoexistsenceHighProb}
	as well as in Proposition \ref{prop:CoexistencewithZeroParticle}, i.e., 
	initially each dormant site has sufficiently many particles with high probability,   
	we have equality of the limiting shapes for both the types 
	along the directions of the diagonals $( \pm 1/d, \cdots, \pm 1/d)$. 
	For discrete time frog model, the limiting shape 
	is always contained in the full diamond $B_1(\mathbf{0}, 1) := 
	\{ \x \in \R^d : ||\x||_1 \leq 1 \}$, the $L_1$ unit ball.   
	For supercritical oriented percolation, almost surely there are infinite oriented open paths along the directions of diagonals (Theorem 1.3 of \cite{R02}). 
	Our coupling with oriented percolation ensures that along those infinite oriented open paths, active particles take oriented steps without being lazy and hence, the limiting shape must coincide with $D$ along the diagonals. Hence, in the set up of Theorem \ref{thm:CoexistsenceHighProb} and 
	Proposition \ref{prop:CoexistencewithZeroParticle},  the limiting shapes 
	for frog models with laziness parameters $p_1$ and $p_2$ coincide with the full diamond 
	$B_1(\mathbf{0}, 1)$ along the diagonals. We don't think that in this set up, 
	we have equality of the complete limiting shapes for $p_1 \neq p_2$. 
	But we don't have a proof at the moment.
\end{remark}

Motivated by the question of coexistence of two types of particles, one can ask the question of coexistence for more than $2$ types of particles. In the next section we explore such questions.

\section{ Coexistence of $l \geq 2$ many types}
\label{sec:Coex_Multi-types}

Main results of this section show coexistence of competing discrete time frog models
for more than two types. The first result shows coexistence of $2^d$ many types of frogs on $\Z^d$ for $d\geqslant 2$. In particular, we can show coexistence of $4$ types of particles on $\Z^2$. 
Further, we show an instance of coexistence of infinitely many types of frogs on $\Z^d$ for $d \geq 3$. All these results are applicable for general laziness parameters provided  
each dormant site has sufficiently large number of particles with high probability.

%We start with an initial configuration of i.i.d. non-negative integer valued random variables $\{\eta(\x) : \x \in \Z^d\}$ independent of the collection $\{(I^{\x,j}_n, U^{\x,j}_n) : \x \in \Z^d, j \geq 1\}_{n \geq 1}$. Now we state the main results of this section about multi-type coexistence. 

\begin{theorem}
	\label{thm:CoexistencewithMultipleTypes}
	\begin{itemize}
		\item[(i)] On $\Z^d$ for $d \geq 2$, we start with $2^d$ types. For $1 \leq i \leq 2^d$ the laziness parameter corresponding to the $i$-th type is given by $p_i \in (0,1]$. Let $\tilde{p} = \min\{p_1, \cdots, p_{2^d}\}$. Fix $2^d$ many distinct points $\x_1, \cdots, \x_{2^d} \in \Z^d$ such that  for all $1 \leq i \leq 2^d$, the site $\x_i$ has a single active type $i$ particle at time zero. The corresponding conditional probbaility measure is denoted as $\P_{\x_1, \cdots, \x_{2^d}}(\cdot)$. There exist $M = M(\tilde{p}, d) \in \N$  and $ \gamma_1 \in (0,1)$ such that if $\P(\eta(\mathbf{0})\geq M) \geq \gamma_1$, then coexistence 
		probability (of all $2^d$ types) is positive w.r.t. $\P_{\x_1, \cdots, \x_{2^d}}(\cdot)$.
		
		\item[(ii)] {On $\Z^d$ for $d \geq 3$}, we consider infinitely many 
		types with laziness parameters $p_i $ for $i \geq 1$ such that $\inf\{p_i : i \geq 1\} = p_0 \in (0,1]$. We start from infinitely many distinct sites $\x_i: i \geqslant 1$ each having an active particle of type $i$. We further require a condition that for all $i_1 , i_2\geq 1$ we have $\x_{i_1}(j) = \x_{i_2}(j)$ for all $1 \leq j \leq d-1$. 
		Then there exists $M = M(p_0, d) \in \N$ and $ \gamma_2 \in (0,1)$ such that if $\P(\eta(\mathbf{0})\geq M) \geq \gamma_2$, then  there is a positive probability of coexistence of infinitely many types.   
	\end{itemize}
\end{theorem}

\noindent \textbf{Proof:} For (i) the idea of the proof is very similar to that of Theorem \ref{thm:CoexistsenceHighProb} or Proposition \ref{prop:CoexistencewithZeroParticle}. 
We present only a sketch here. 
We choose $m$ such that $||\x_i||_1  \leq m$ for all $1 \leq i \leq 2^d$. 
We consider a bijective map 
$$
f: \{ \x_i : 1 \leq i \leq 2^d\} \mapsto \{ \theta : \theta \in \{ + 1, - 1 \}^d \}. 
$$
Set $k_i := ||\x_i - mf(\x_i)||_1$ for $1 \leq i \leq 2^d$. Following the notation in Section \ref{sec:CouplingOP}, we consider the $f(x_i)$ oriented orthants $\Lambda^{f(\x_i)}(m f(\x_i))$ for $1 \leq i \leq 2^d$, and by an abuse of notation we shall call them $\Lambda_i$ for $1 \leq i \leq 2^d$ respectively. 
For each $1 \leq i \leq 2^d$, we  consider a path 
$\pi_i := \{\x_{i0} := \x_{i},\x_{i1}, \cdots, \x_{ik_i} = mf(\x_i) \}$ of length 
$k_i$ from $\x_i$ to $mf(\x_i)$ such that $||\x_{ij}||_1 < m$ for all $0 \leq j \leq k_i - 1$. 
Choose $M = M(\tilde{p}, d) \in \N$ and $\gamma_1 \in (0,1)$ so that $
\gamma_1 g(M, \tilde{p}) > p^\uparrow_c$.
Following the coupling in Section \ref{sec:CouplingOP}, for the orthant $\Lambda_i$
we consider the oriented percolation model for the choice $p = \tilde{p}, M = M(\tilde{p}, d)$ and $\ell = k_i$.

We consider the event that for all $1 \leq i \leq 2^d$, 
the initially active particle of $i$-th type follows the path 
$\pi_i$ and reaches the site $mf(\x_i)$ exactly at time $k_i$. 
We further require that for all $1 \leq i \leq 2^d$, the site $mf(\x_i)$ is a percolating 
point, i.e., it has an infinite oriented open path in the orthant $\Lambda_i$.
We need to control the movement of the newly activated particles also. 
Let $L_n(i)$ denote the set of location(s) of all $i$-type particles at time $n$. 
The event $\cap_{i=1}^{2^d} \cap_{n=0}^{k_i} ( L_n(i) = \{\x_{in}\})$ 
controls
movement of all the newly activated particles. 
We observe that for any $1 \leq i \neq j \leq 2^d$ the orthants $\Lambda_i$
and $\Lambda_j$ are disjoint. 
Finally, the same argument as in Theorem \ref{thm:CoexistsenceHighProb} 
gives us that intersection of all the above events 
is of positive probability and implies coexistence of all $2^d$ types.       

%all the sites in the set $B_\infty(\mathbf{0}, m)$, which become activated due  to the movement of initially active particles, stay still till time $2md + 1$ and for each $\theta \in \{1,2\}^d$ the point $\theta_m  + (-1)^{\theta(1)}e_1$ is a percolating point, i.e., $\theta_m  + (-1)^{\theta(1)}e_1$
%has an infinite oriented open path in $\Lambda^{\theta_m}$. 
%The same argument as in Theorem \ref{thm:CoexistsenceHighProb} gives us that the event described above 
%would imply coexistence of all the $2^d$ types of particles.  
%Further, we can choose $M \in \N$ and $\gamma_1 \in (0,1)$ so that coexistence occurs with positive probability. 

For (ii), our assumption ensures that  we must have 
$\x_{i_1}(d) \neq \x_{i_2}(d)$ for all $i_1 \neq i_2$.

Let us recall that we start with one active particle of type $i$ from the site 
$\x_i$ at time $0$, and let $\P_{\bigotimes_i \x_i}(.)$ 
denote the corresponding conditional probability measure. 
For $i \geq 1$ we consider the oriented orthant 
$$
\Lambda^+_i := \{ \y \in \Z^d : \y(d) = \x_i(d), 
\y(j)\geq \x_i(j) \text{ for all }1 \leq j \leq d - 1\}.
$$  
We can chose $M = M(p_0, d)\in \N$, where $p_0 = \inf\{ p_i : i \geq 1\} > 0$, large so that   
$$
g(M, p_0) > p_c^{\uparrow}(d-1)
$$ where $p_c^{\uparrow}(d-1)$ is the  threshold for 
oriented site percolation on $\mathbb{Z}^{d-1}$.
Clearly, we can chose $\gamma_2 \in (0,1)$ such that $\gamma_2 g(M, p_0) > p_c^{\uparrow}(d-1)$. 
Now we refer to our coupling construction in Section \ref{sec:CouplingOP}.
For each $i \geq 1$ for the orthant $\Lambda^+_i$ we consider
the oriented percolation model with $p = p_0, 
M = M(p_0, d)$  and $\ell   = 0$. 

For $i \geq 1$ we define the event $E_i$ as 
\begin{align*}
E_i & := \bigl \{ \text{The point } \x_i  
\text{ has an infinite open path in the set } \Lambda^+_i \bigr \}.
\end{align*}
We observe that for $i, i^\prime \geq 1$ with $i \neq i^\prime$ events $E_i$ and $E_{i^\prime}$
are supported on disjoint set of random vectors implying that the events 
$E_l : l \geq 1$ are mutually independent. 
Clearly, we have $\P_{\bigotimes_i \x_i} (E_1)=\P_{\bigotimes_i \x_i} (E_l)$ for all $l \geq 1$
and our choice of $\gamma_2, M$ ensures that $\P_{\bigotimes_i \x_i} (E_1)> 0$. 
Hence, we have $\sum_{l=1}^{\infty} \P_{\bigotimes_i \x_i} (E_l) =\infty $. 
An application of second Borel Cantelli lemma gives us that 
$E_l$'s happen infinitely often with probability $1$. 

Now, Same argument as in Theorem \ref{thm:CoexistsenceHighProb} gives us that
on the event $E_i$, any (open) site on an infinite oriented open path
starting from $\x_i$ and on the set $\Lambda^+_i \cap \{ \x \in \mathbb{Z}^d : \x(d) = \x_i(d) \}$
must be activated by an $i$ type particle only. 
This implies that, for $d\geq 3$ infinite number of types can co-exist. 
\qed

\end{document}